\documentclass[10pt]{amsart}

\usepackage{hyperref}
\usepackage{xcolor}
\usepackage{amssymb}
\usepackage{amsmath}

\let\newpf\proof \let\proof\relax

\def\bm{\begin{matrix}}
\def\em{\end{matrix}}

\newcommand{\bt}{\begin{thm}}
\newcommand{\et}{\end{thm}}

\newcommand{\bl}{\begin{lemma}}
\newcommand{\el}{\end{lemma}}

\newcommand{\beq}{\begin{eqnarray}}
\newcommand{\eeq}{\end{eqnarray}}

\def\be{\begin{equation}}
\def\ee{\end{equation}}

\def\ba{{\begin{align}}}
\def\ea{{\end{align}}}

\def\0{{\mathbf 0}}

\def\cal{\mathcal}

\newtheorem{thm}{Theorem}[section]

\newtheorem{lem}[thm]{Lemma}
\newtheorem{lemma}[thm]{Lemma}

\theoremstyle{remark}
\newtheorem{rmk}{Remark}[section]

\newtheorem{example}{Example}[section]

\numberwithin{equation}{section}

\def \bn {\hfill \\ \smallskip\noindent}

\theoremstyle{definition}

\def\proof{\bn {\bf Proof.} }

\def\note#1
{\marginpar

{\tiny $\leftarrow$
\par
\hfuzz=20pt \hbadness=9000 \hyphenpenalty=-100 \exhyphenpenalty=-100
\pretolerance=-1 \tolerance=9999 \doublehyphendemerits=-100000
\finalhyphendemerits=-100000 \baselineskip=6pt
#1}\hfuzz=1pt}

\newcommand{\C}{{\mathbb C}}

\newcommand{\M}{{\mathbb M}}
\newcommand{\N}{{\mathbb N}}
\newcommand{\Q}{{\mathbb Q}}
\newcommand{\R}{{\mathbb R}}
\newcommand{\T}{{\mathbb T}}

\newcommand{\Z}{{\mathbb Z}}

\def\B0{{\bold{0}}}

\catcode`\@=12

\def\Empty{}
\newcommand\oplabel[1]{
  \def\OpArg{#1} \ifx \OpArg\Empty {} \else
  	\label{#1}
  \fi}

\newcommand{\comm}[1]{}
\newcommand{\comment}[1]{}

\begin{document}

\title{Singular continuous spectrum and generic full spectral/packing dimension for unbounded quasiperiodic Schr\"odinger operators}

\author{Fan Yang and Shiwen Zhang}

%\address{}

%\email{a}

%\address{a}
%\email{a}
%\thanks{S.J. is a 2014-15 Simons Fellow. This research was partially supported by  and .}

\begin{abstract}

We proved that Schr\"odinger operators with unbounded potentials
$(H_{\alpha,\theta}u)_n=u_{n+1}+u_{n-1}+
\frac{g(\theta+n\alpha)}{f(\theta+n\alpha)} u_n$ have purely singular
continuous spectrum on the set $\{E: 0<L(E)<\delta{(\alpha,\theta;f,g)}\}$,
where $\delta$ is an explicit function and $L$ is the Lyapunov
exponent. We only require $f,g$ are H\"older continuous functions and $f$ has finitely many zeros with weak non-degenerate assumptions. Moreover, we show that for generic $\alpha$ and a.e. $\theta$, the spectral measure of $H_{\alpha,\theta}$ has full spectral/packing dimension. %This extends results in \cite{JY} for meromorphic potentials and work of \cite{JZ} to unbouned Schr\"odinger operators. 

\end{abstract}

\maketitle

\section{Introduction and main results}
Unbounded Schr\"odinger operators attract a lot of interests in both physics and math literatures, see e.g. \cite{GFP,BLS,B,FP,S85,SS,KMP90,KMP,GJMS,JNS,GKDS,JLiu,JY}. In general, high barriers of the operators lead to the absence of absolutely continuous spectrum or even pure point spectrum with localized eigenstates, see e.g. \cite{SS,KMP,GJMS,JNS}. One significant example of unbounded quasiperiodic Schr\"odinger operator is the Maryland model, proposed by Grempel, Fishman, and Prange \cite{GFP} as a linear version of the quantum kicked rotor. Mathematically, it is interesting due to its richness of spectral theory, see Simon \cite{CFKS,S85}. A complete description of spectral transitions with respect to all parameters for the Maryland model was given recently by Jitomirskaya and Liu in \cite{JLiu}. There are also a number of multidimensional generalizations of Maryland model, see e.g. \cite{BLS,FP}. % and other related models, for example, surface Maryland model (e.g. [20]).
On the other hand, there is a simple way of excluding the existence of point spectrum that relies on Gordon condition (a single/double almost repetition) of the potential going back to Gordon \cite{G76}; see also in \cite{CFKS,D00}. Most applications of the Gordon-type argument quantify the competition between the quality of repetitions and the Lyapunov growth, and are concentrated ` on bounded potentials, see e.g. \cite{D00,D04,BD}. Simon   in \cite{S85} applied Gordon-type of arguments to Maryland model and obtainded purely singular continuous spectrum for generic frequencies. Jitomirskaya and Liu in \cite{JLiu} obtained the sharp parameters region (both in phase and frequency) for singular continuous spectrum by the refined Gordon-type argument (see Theorem \ref{thm:Gordon}).
%On the other hand, for Schr\"odinger operators with unbounded potential, Simon and Spencer in \cite{SS} established a simple general criterion for the absence of an absolutely continuous spectrum. Later on, Kirsch, Molchanov and Pastur in \cite{KMP} gave a finer criterion on the unboundedness to obtain pure point spectrum of the half-line Schr\"odinger operator for a.e. bounary condition. The result in \cite{SS} says that so long as there are arbitrarily high barriers, the Schr\"odinger operator  has no absolutely continuous spectrum. There are also strong interests of studying operators with high barriers, i.e., unbounded potentials.   It is natrual to ask the presence of singular continuous spectrum for unbounded potential.  It is an exactly solvable example, e.g. [7]. More recently, in the context of one-dimentional quasiperiodic Sch\"odinger operators, Jitomirskaya-Liu \cite{JLiu} and
This sharp result about the singular continuous spectrum was generalized by Jitomirskaya and Yang in \cite{JY} to unbounded quasiperiodic Schr\"odinger operators with meromorphic potentials.  
%In particular, the sharp phase-frequency transition line between the pure point spectrum and purely singular continuous spectrum was obtianed for the Maryland model in \cite{JLiu}. 
%The key ingredient to apply Gordon-type argument to unbounded potentials is a Lyapunov growth estimate for unbounded linear cocycles. 
%The meromorphic potentials (analycity of the singularities) play important role in both work. 

Recently, there has been an increased interest in developing methods that don't involve analyticity, see e.g. \cite{K,B07,WZ,JK}. The methods to exclude point spectrum, used in \cite{S85,JLiu,JY}, strongly rely on the meromorphic potentials, where the singularities are analytic. It is still not clear how to obtain singular continuous spectrum in the best possible arithmetic regime for rough unbounded potential. This leads to our first motivation. In this paper, we study unbounded quasiperioidc Schr\"odinger operators with lower regularity assumptions. %We prove a weakening uniformity of upper semicontinuity of the Lyapunov exponent for unbounded cocyles. 
We obtain the absence of point spectrum in the sharp parameters regime. This generalizes the results of \cite{JY}  to more general unbounded potentials. 

 Moreover, we want to study finer decompositions of the singular continuous measure of quasiperiodic Schr\"odinger operators with rough unbounded potentials. In the recent work of Jitomirskaya and Zhang in \cite{JZ},  a quantitative version of Gordon-type results was obtained: a quantitative strengthening  (multiple almost repetitions) implies quantitative continuity of the spectral measure. While the main application to the quasiperiodic case also requires the potential to be bounded. %In the spirit of the key lemma that we used to show the absence of point spectrum,
 For unbounded operators, the difficulty lies in obtaining uniform upper semicontinuity of the Lyapunov exponent. This is resolved in Lemma \ref{lem:infA}, which leads to quantitative spectral continuity after incorporating the abstract scheme developed in \cite{JZ}. The main consequence is the full spectral/packing dimension of the spectral measures for generic frequencies, which generalizes the results of \cite{JZ} to unbounded potentials.   

More precisely, we study lattice Schr\"odinger operators on ${\ell}^2(\Z)$ of the form:
\begin{equation}\label{eq:op}
(H_{\alpha,\theta}u)_n=u_{n+1}+u_{n-1}+\frac{g(\theta+n\alpha)}{f(\theta+n\alpha)} u_n,
\end{equation}
where $\alpha \in \R\setminus \Q$ is the
frequency, $\theta \in \T:=\R/\Z$ is the phase.  Let $C^{\tau_0}(\T,\R)$ be the space of $\tau_0$-H\"older continuous functions with norm:
\begin{align}\label{def:norm}
\|f\|_{\tau_0}:=\sup_{x\in\T}|f(x)|+\sup_{x,y\in\T}\frac{|f(x)-f(y)|}{|x-y|^{\tau_0}}<\infty
\end{align}
for some $0<\tau_0\le 1$. 
We assume $f,g\in C^{\tau_0}(\T,\R)$. 
We also require $f$ to be in the following space of functions introduced in \cite{HYZ}:
\begin{align}\label{def:F}
{\cal {F}}(\T,{\R}):=\Big\{&f\in C^0({\T,{\R}}):\ \exists m \in\N^+, \theta_\ell\in\T, 0<\tau_\ell\le 1,\ell=1,\cdots,m \nonumber
\\
&\qquad \textrm{ such that}\ \ 
\widetilde f(\theta):=\frac{f(\theta)}{\prod_{\ell=1}^m|\sin \pi(\theta-\theta_{\ell})|^{\tau_{\ell}}}
\in C^0(\T,{\R})
\ {\rm and}\ \inf_\T |\widetilde f(\theta)|>0. 
\Big\} 
\end{align}

Let $\frac{p_n}{q_n}$ be the continued fraction approximants of $\alpha\in \R \setminus \Q$. We define index $\delta$ as follows:

 \begin{equation} \label{def:delta}
  \delta( \alpha,\theta)= \delta( \alpha,\theta;f,g) = \limsup_{n\rightarrow \infty} \dfrac{\sum_{\ell=1}^m \tau_{\ell}\ln\|q_n(\theta-\theta_i) \|_{\R / \Z}+ \tau_{\rm min} \ln q_{n+1}} {q_n},
 \end{equation}   
where $\tau_{\rm min}=\min_{0\le\ell\le m}\tau_\ell$ and 
$\|x\|_{\R / \Z}=\min_{\ell\in \Z}|x-\ell|$.  Let $L(E)$ be the Lyapunov
exponent, see (\ref{def:lyp}). $L(E)$ depends also on $\alpha$ but we suppress
it from the notation as we keep $\alpha$ fixed. 

Our main result is:
 \begin{thm}\label{thm:sc}
 	Let $f,g$ be $\tau_0$-H\"older continuous, $f$ be given as in \eqref{def:F} and let  $ \delta(\alpha,\theta;f,g) $ be as in (\ref{def:delta}). Then $ H_{\alpha, \theta} $ %(\ref{eq:op}) 
has no eigenvalues on 
{$\{ E:L(E)< \delta(\alpha, \theta)\}$. }
\end{thm}
\begin{rmk}
	 Absence of absolutely continuous spectrum follows
 from a.e. positivity of the Lyapunov exponents and holds
 for all unbounded potentials \cite{SS}. Then  $ H_{\alpha, \theta} $  has purely singular continuous spectrum on
 $\{ E:0<L(E)< \delta(\alpha, \theta)\}$.
\end{rmk}

It is easy to check that the class ${\cal {F}}(\T,{\R})$ in \eqref{def:F} contains all $C^k,k\ge1$ continuous potentials with finitely many non-degenerate zeros
\footnote{We say $\theta_0\in\T$ is a non-degenerate zero of $f\in C^{k}(\T,\C)$  if $f(\theta_0)=0$ and $f^{(k)}(\theta_0)\neq0 $.}. In particular, it covers the potential studied in \cite{JY}, where $g$ is assumed to be Lipschitz continuous and $f$ is analytic, where all $\tau_\ell=1, \ell=0,\cdots, m$.  The following examples where the potential is only H\"older continuous show that such generalization is indeed nontrivial.

\begin{example}
	Let $
	(H_{\alpha,\theta}u)_n=u_{n+1}+u_{n-1}+\frac{1}{f(\theta+n\alpha\, {\rm mod}\,1)} u_n,
$
	where 
	\begin{align}
	f(\theta)=\begin{cases}
0 &,\ \ \theta=0; \\
\theta\cos(\frac{1}{\theta})+2\sqrt{\theta} &, \ \ 0<\theta\le\frac{1}{2};\\
f(1-\theta) &, \ \ \frac{1}{2}<\theta\le 1.
	\end{cases}
	\end{align}
It is an easy exercise by definition that $\theta \cos(\frac{1}{\theta})$	is sharp $\frac{1}{2}$-H\"older continuous near $0$. Moreover, $\frac{f(\theta)}{\sqrt{|\sin\pi\theta|}}\ge 1$ and is continuous on $\T$. Therefore, $f\in{\cal F}$ with $\tau_1=\frac{1}{2}$.  
	Let \begin{align}\label{def:tildedelta}
	\widetilde\delta( \alpha,\theta) = \limsup_{n\rightarrow \infty} \dfrac{\ln\|q_n\,\theta \|_{\R / \Z}+  \ln q_{n+1}} {{q_n}}.
	\end{align}
Theorem \ref{thm:sc} implies that  $ H_{\alpha, \theta} $  has purely singular continuous spectrum on
$\{ E:0<L(E)< \frac{1}{2}\widetilde\delta(\alpha, \theta)\}$.
\end{example}  

\begin{example}
	Let 
	\begin{equation}
	(H_{\alpha,\theta}u)_n=u_{n+1}+u_{n-1}+{h(\theta+n\alpha \, {\rm mod}\,1)} u_n,
	\end{equation}	
	where 
	\begin{align}
	h(\theta)=\begin{cases}
	\log\theta &,\ \ 0<\theta\le\frac{1}{2};\\
\log(1-\theta) &, \ \ \frac{1}{2}<\theta< 1.
	\end{cases}
	\end{align}
Let $	\widetilde \delta( \alpha,\theta) $ be as in \eqref{def:tildedelta}. 
Let \begin{align}
g(\theta)=\begin{cases}
0 &,\ \ \theta=0; \\
\theta\log{\theta} &, \  \ 0<\theta\le\frac{1}{2};\\
g(1-\theta) &, \ \ \frac{1}{2}<\theta\le 1.
\end{cases}
\end{align} Clearly, $h(\theta)=\frac{g(\theta)}{\theta}$ on $(0,1)$ and $g(\theta)$ is $(1-\epsilon)$-H\"older continuous on $[0,1]$ for any $0<\epsilon<1$. By Theorem \ref{thm:sc},  $ H_{\alpha, \theta} $  has purely singular continuous spectrum on
	$\{ E:0<L(E)< (1-\epsilon)\widetilde\delta(\alpha, \theta)\}$ for any $0<\epsilon<1$. It is also easy to check that exact the same thing holds true for $\widetilde h(\theta)=\chi_{(0,1/2]}\log\theta-\chi_{(1/2,1)}\log(1-\theta)$. 
\end{example}  

 Once we obtain singular continuous spectral measure, we want to move further to study the factal dimension properties of the measure. Let  $\mu=\mu_{\alpha,\theta}$  be the spectral measure of the operator \eqref{eq:op}. The fractal properties of $\mu$ are closely related to the  
boundary behavior of its Borel transforms $M(E+ i\varepsilon)=\int\frac{{\rm d}\mu(E')}{E'-(E+ i\varepsilon)}$, see e.g. \cite{DJLS}.  We are interested in the following local fractal dimension of $\mu$. For any compact set $I\subset \R$, let $\mu_I$ be the restriction of $\mu$ on $I$.	We say $\mu_I$ is $\gamma$-spectral continuous if for some $\gamma\in(0,1)$ and $\mu$ a.e. $E\in I$, we have
	\begin{equation}\label{def:speconti}
	\liminf_{\varepsilon\downarrow0}\varepsilon^{1-\gamma}|M(E+ i\varepsilon)|<\infty.
	\end{equation}
	Define the (upper) spectral dimension  of $\mu_I$ to be
	\begin{align}\label{def:dimspe}
	{\rm dim}_{\rm spe}(\mu_I)=\sup\big\{\gamma\in(0,1):\ \mu_I\ \textrm{is}\ \gamma\textrm{-spectral continuous}\big\}.
	\end{align}
Consider another arithmetic index of $\alpha$  defined by $\beta(\alpha)=\limsup_{n\rightarrow \infty} \dfrac{ \ln q_{n+1}} {q_n}.$ It is well known that for a.e. $\theta\in\T$, $ \delta( \alpha,\theta)=\tau_{\rm min}\cdot \beta(\alpha)$.
 Combine the Lyapunov growth estimates (see Lemma \ref{lem:infA}) and the recent work of Jitomirskaya-Zhang \cite{JZ}, we have the following quantitative spectral continuity for $\mu_I$.  
\begin{thm}\label{thm:dim}
Under the same assumption of Theorem \ref{thm:sc}. For any compact set $I\subset \R$, there is a constant $C=C(\|f\|_\infty,\|g\|_\infty,I)$ such that for a.e. $\theta\in\T$ if $\tau_{\rm min}\cdot \beta(\alpha)>C$, then $\mu_I$ has positive spectral dimension. In particular, for a.e. $\theta$, if $\beta(\alpha)=\infty$, then $\mu_I$ has full spectral  dimension. 
\end{thm}
\begin{rmk}
The Hausdorff/packing dimension of a (Borel) measure $\mu$, namely, ${\rm dim}_{\rm H}(\mu)/{\rm dim}_{\rm P}(\mu)$ is defined through the  $\limsup/\liminf$ ($\mu$ almost everywhere) of its $\gamma$-derivative $\lim_{\varepsilon\downarrow 0}\frac{\ln \mu(E-\varepsilon,E+\varepsilon)}{\ln \varepsilon}.$ It is well known (see e.g. \cite{Fal,JZ}) that the (local) packing dimension is bounded from below by the spectral dimension in \eqref{def:dimspe}. Therefore, we have that for $\beta(\alpha)=\infty$ (which is a generic subset of $[0,1]$), the restriction of $\mu$ on any compact set has full packing dimension. On the other hand,  it was showed in \cite{S07} that $\mu$ is Haurdorff singular (${\rm dim}_{\rm H}(\mu)=0$) on $\{ E:L(E)>0\}$ for any $\alpha,\theta$. Therefore, we have a generic subset of frequencies $\alpha$,  such that for a.e. $\theta$, $\mu_I$ is not exact dimensional (the Haurdorff dimension and the packing dimension do not equal). 
\end{rmk}
 
\section{Preliminaries: cocycle, Lyapunov exponent} 
From now on, we assume $f\in C^{\tau_0}$ and has the expresssion as in \eqref{def:F}:
$$f(\theta)=\widetilde f(\theta)\prod_{\ell=1}^m|\sin \pi(\theta-\theta_{\ell})|^{\tau_{\ell}},\ \ \theta_\ell\in\T,\ \tau_\ell\in(0,1],\ \ell=1\cdots,m,\  \widetilde f\in C^0
\ {\rm and}\ \inf_\T |\widetilde f(\theta)|>0. $$  
Let $S=\bigcup_{\ell=1}^m\theta_\ell+\alpha\Z+\Z$. We fix $E$ in the spectrum and $\theta\in S^c$.

A formal solution of the equation  $ H_{\alpha, \theta}u=Eu $ can be reconstructed via the following equation
 \begin{align}\label{def:A}
\left (\begin{matrix}  u_{n+1} \\ u_{n} \end{matrix} \right )= A(\theta+n\alpha)\left (\begin{matrix}  u_{n} \\  u_{n-1}  \end{matrix} \right ),
 \end{align} 
 where
$
A(\theta)=\left (\begin{matrix}  E-\frac{g(\theta)}{f(\theta)}   &  -1 \\ 1  &   0 \end{matrix} \right )
$ is the so-called transfer matrix.

 The pair $(\alpha,A)$ is the cocycle corresponding to the operator  (\ref{eq:op}). 
 It can be viewed as a linear skew-product $(x,\omega)\mapsto(x+\alpha,A(x)\cdot \omega)$.
 Generally, one can define $M_n$ for an invertible cocycle $(\alpha, \M)$ by $(\alpha,M)^n=(n\alpha,M_n)$, $n \in Z$ so that for $n \geq 0$:
 \begin{align}\label{def:Mn}
 M_n(x)=M(x+(n-1)\alpha)M(x+(n-2)\alpha) \cdots M(x),
 \end{align}
and $M_{-n}(x)=M^{-1}_n(x-n\alpha)$. 

The Lyapunov exponent of a cocycle $(\alpha,M)$ is defined by$$L(\alpha,M)=\lim_{n \rightarrow \infty}\frac{1}{n} \int_{\T} \ln \|M_n(x)\|\mathrm{d}x.$$

Let $A(x)=\frac{1}{f(x)}D(x)$, where 
\begin{align} \label{def:AD}
D(x)=
\left(
\begin{matrix}
    E f(x)- g(x)        &  -f(x)
 \\       f(x)              &   0
 \end{matrix}
\right)
 \end{align} 
is the regular part of $A(x)$.
It is convenient to assume  $\int_{\T} \ln{|f(x)|} \mathrm{d}x=0$, otherwise it is enough to consider the constant scaling by $b=e^{-\int_{\T} \ln{|f(x)|} \mathrm{d}x}$ and $f_b:=bf(\theta),g_b:=bg(\theta)$ in the potential. We have 
\begin{equation}\label{def:lyp}
L(E):=L({\alpha, A} )= L(\alpha, D).
\end{equation}

  \begin{lemma}\label{lem:AJ09}\cite{AJ1}
 Let $\alpha\in \R\backslash\Q $,\ $\theta\in\R$ and $0\leq j_0 \leq q_{n}-1$ be such that 
 $$\mid \sin \pi(\theta+j_{0}\alpha)\mid = \inf_{0\leq j \leq q_{n}-1} \mid \sin \pi(\theta+j\alpha)\mid ,$$
 then for some absolute constant $C>0$,
 $$-C\ln q_{n} \leq \sum_{j=0,j\neq j_0}^{q_{n}-1} \ln \mid \sin \pi (\theta+j\alpha) \mid+(q_{n}-1)\ln2 \leq C\ln q_n.$$
 \end{lemma}
 
 We will also use that the denominators of continued fraction approximants of $\alpha$ satisfy
 
 $$ \| k\alpha \|_{\R \backslash \Z} \geq \| q_n \alpha \|_{\R \backslash \Z}, 1 \leq k < q_{n+1}, $$ 
 and
\begin{equation}\label{contfraction}
 \dfrac{1}{2q_{n+1}}  \leq  \|q_n \alpha\|_{\R \backslash \Z}  \leq  \dfrac{1}{q_{n+1}}.
\end{equation} 

A quick corollary of subadditivity and unique ergodicity is the 
following uniform upper semicontinuity result:
 \begin{lemma}[\cite{F}, see also in \cite{AJ2,JM}]\label{lem:F} 
Suppose $(\alpha,A)$ is a continuous cocycle. Then for any $ \varepsilon>0$, there exists $C(\varepsilon)>0$, such that for any $x\in \T$ we have
$$\|A_n(x) \| \leq C e ^{n(L(A)+\varepsilon)}.$$
\end{lemma}
\begin{rmk}
	The above result was first obtained by Furman in \cite{F} in a relatively more general setting. Jitomirskaya and Mavi generalzied the result to the case where $A(x)$ is bounded in norm (from above) and has points of discontinuity with Lebesgue measure zero. In general, such uniform ( in $x$) $\limsup$ may not hold if $A(x)$ is unbounded. In section \ref{sec:dim}, we obtain a weakened version (see Lemma \ref{lem:infA}) of such upper semicontinuity for quasiperiodic unbounded cocycles. The estimate plays important role in the proof of Theorem \ref{thm:dim}. It is also  of independent interest in the study of the multiplicative ergodic theorem for uniquely ergodic systems. 
\end{rmk}
\begin{rmk}\label{rmk:1-co}
Applying this to 1-dimensional continuous cocycles, we get that if $g$
is a continuous function such that $\ln |g| \in L^1(\T)$,
then $$|\prod_{j=a}^b g (x+j\alpha) | \leq e^{(b-a+1)(\int \ln|g| {\rm d}\theta+\varepsilon)}.$$
\end{rmk}

\section{Absence of Point Spectrum}
By the definition of $\delta(\alpha, \theta)$, for any $\varepsilon>0 $, there exists a subsequence $q_{n_i}$ of $q_n $ such that
\begin{align}\label{eq:q_n}
\dfrac{\sum_{\ell=1}^m \tau_\ell\ln\|q_{n_i}(\theta-{\theta_\ell}) \|+ \tau_{\rm min}\ln q_{n_i+1}} {q_{n_i}} >\delta(\alpha,\theta)-\frac{\varepsilon}{4}.
\end{align}

In this section, we omit the subindex $n_i$ and still denote the subsequence satisfying (\ref{eq:q_n}) by $q_n$ whenever it is clear.
\subsection{key lemmas}
The following lemmas are in spirit close to Lemma 4.2 in \cite{JY}. 
For all $1\le \ell  \le m$, 
let $j_{{\ell}}\in [0,q_n)$ be such that the following holds:
$|\sin\pi\big(\theta-\theta_{\ell}+j_{{\ell}} \alpha \big)| 
=\inf_{0\leq j < q_n}|\sin\pi\big(\theta-{\theta_\ell}+j\alpha\big)|.$

\begin{lemma} \label{lem:key}
	
	If $\delta(\alpha,\theta) > 0$, then 
	\begin{align}\label{eq:key-1}
	\prod_{\ell=1}^m | \sin \pi (\theta-\theta_\ell+j_{{\ell}} \alpha) |^{\tau_{\ell}} \geq \frac{e^{q_{n}(\delta-\frac{\varepsilon}{2})}}{q^{\tau_{\rm min}}_{n+1}},
	\end{align}
	where $\tau_{\rm min}=\min_{0\le\ell\le m}\tau_\ell$.
\end{lemma} 
\proof
	By (\ref{eq:q_n}), 
	\begin{align}
	\prod_{\ell=1}^m\|q_{n}(\theta-\theta_\ell)\|^{\tau_\ell} > \frac{e^{q_{n} (\delta-\frac{\varepsilon}{4})}}{q^{\tau_{\rm min}}_{n+1}}.
	\end{align}
	
	In particular, 
	\begin{align}
	\|q_{n} (\theta -{\theta_\ell})\|^{\tau_\ell} >  \frac{e^{q_{n} (\delta-\frac{\varepsilon}{4})}}{q^{\tau_{\rm min}}_{n+1}}>  \frac{e^{q_{n} (\delta-\frac{\varepsilon}{4})}}{q^{\tau_{\ell}}_{n+1}}
	\end{align}
  for any $1\leq \ell \leq m$.
	Therefore, $\|q_{n} (\theta -\theta_\ell)\|\ge q_{n+1}^{-1} e^{\tau_\ell^{-1}q_{n} (\delta-\frac{\varepsilon}{4})}\ge 2q_nq_{n+1}^{-1} $ and 
	
	\begin{align*}
	q_n| \sin \pi (\theta-{\theta_\ell}+j_{{\ell}} \alpha)|
	\geq 2q_n\|\theta-{\theta_\ell}+j_{{\ell}} \alpha\| 
	\geq 2\|q_{n}(\theta-{\theta_\ell}+j_{{\ell}} \alpha)\| 
	\geq &{2\|q_{n}(\theta-{\theta_\ell})\|-2j_{{\ell}}\|q_n\alpha\|} \\
	\ge & {2\|q_{n}(\theta-{\theta_\ell})\|-2q_nq_{n+1}^{-1}} \\
	\geq &\|q_{n}(\theta-{\theta_\ell})\| ,
	\end{align*}
	provided $q_n$ large. Denote 
	\begin{align}\label{eq:tausum}
\tau_{\rm sum}=\sum_{\ell=1}^m\tau_{\ell}.
	\end{align}
	Then we have 
\begin{align*}
\prod_{\ell=1}^m | \sin \pi (\theta-\theta_\ell+j_{{\ell}} \alpha) |^{\tau_{\ell}}
\ge \prod_{\ell=1}^m  \frac{\|q_{n}(\theta-{\theta_\ell})\|^{\tau_{\ell}}}{q^{\tau_{\ell}}_{n}}
&
= \prod_{\ell=1}^m  \|q_{n}(\theta-{\theta_\ell})\|^{\tau_{\ell}}\,\prod_{\ell=1}^m  \frac{1}{q^{\tau_{\ell}}_{n}}\\
&	> \frac{e^{q_{n} (\delta-\frac{\varepsilon}{4})}}{q^{\tau_{\rm min}}_{n+1}} \cdot \frac{1}{q_{n}^{\tau_{\rm sum}}}\\
&
> \frac{e^{q_{n}(\delta-\frac{\varepsilon}{2})}}{q^{\tau_{\rm min}}_{n+1}},
\end{align*}
	provided $q_n$ large. 
\qed

\begin{lemma}\label{lem:key2}
	Let $f\in{\cal {F}}$ as defined in (\ref{def:F}) where 
	\begin{align*}
	f(\theta)=\widetilde f(\theta){\prod_{\ell=1}^m|\sin \pi(\theta-\theta_{\ell})|^{\tau_{\ell}}},\ \widetilde f\in C^0,\ \ \inf_\T|\widetilde f|>0.
	\end{align*} 
%	We have for $|k|q_n\le e^{\varepsilon q_n}$, 
%	\begin{align}
%	\prod_{j=kq_n}^{(k+1)q_{n}-1}|f(\theta+j\alpha )| \geq  e^{q_n\int_\T\ln | f(\theta)|{\rm d}\theta}\cdot \%dfrac{e^{q_n(\delta-\varepsilon) }}{q^{\tau_{\rm min}}_{n+1}},
%	\end{align}
 Assume further that 
	 $\int_\T\ln | f(\theta)|{\rm d}\theta=0$,  then 
	\begin{align} \label{eq:JY-1}
\prod_{j=0}^{q_{n}-1}|f(\theta+j\alpha )| \geq  \dfrac{e^{q_n(\delta-\varepsilon) }}{q^{\tau_{\rm min}}_{n+1}},
	\end{align}	where $\tau_{\rm min}=\min_{0\le\ell\le m}\tau_\ell$.
\end{lemma}

\proof
By Lemma \ref{lem:AJ09}, 
\begin{align}
\prod_{j=0,j\neq j_{{\ell}}}^{q_n-1}|\sin\pi(\theta-\theta_{\ell}+j\alpha)|
\ge e^{-(q_n-1)\ln2-C\ln q_n} \ge e^{-q_n(\ln2+\frac{\varepsilon}{4}\tau_{\rm sum}^{-1})} \label{eq:AJ09-2},
\end{align}
provided $C\ln q_n<\frac{\varepsilon}{4}\tau_{\rm sum}^{-1} q_n$, in which $C$ is the absolute constant in Lemma \ref{lem:AJ09} and $\tau_{\rm sum}$ is given as in (\ref{eq:tausum}).

Combine \eqref{eq:key-1} with \eqref{eq:AJ09-2}, we have 
	
\begin{align*}
\prod_{j=0}^{q_{n}-1}\prod_{{\ell}=1}^m|\sin\pi(\theta-\theta_{\ell}+j\alpha)|^{\tau_{\ell}}
%=& \prod_{{\ell}=1}^m\,\prod_{j=0}^{q_{n}-1}|\sin\pi(\theta-\theta_{\ell}+j\alpha)|^{\tau_{\ell}}\\
=& \prod_{{\ell}=1}^m\Big(\prod_{j=0,j\neq j_{{\ell}}}^{q_n-1}|\sin\pi(\theta-\theta_{\ell}+j\alpha)|\Big)^{\tau_{\ell}}
\cdot \Big(\prod_{{\ell}=1}^m| \sin \pi (\theta-\theta_\ell+j_{{\ell}} \alpha) |^{\tau_{\ell}}\Big)\\
\ge&  \Big(\prod_{{\ell}=1}^m\big(e^{-q_n(\ln2+\frac{\varepsilon}{4}\tau_{\rm sum}^{-1})} \big)^{\tau_{\ell}}\Big)\,
\cdot\, \frac{e^{q_{n}(\delta-\frac{\varepsilon}{2})}}{q^{\tau_{\rm min}}_{n+1}} \\
=& e^{-q_n(\tau_{\rm sum}\ln2+\frac{\varepsilon}{4})}\,
\cdot\, \frac{e^{q_{n}(\delta-\frac{\varepsilon}{2})}}{q^{\tau_{\rm min}}_{n+1}}\\
=&e^{q_{n}(\delta-\frac{3\varepsilon}{4}-\tau_{\rm sum}\ln2)}\,{q^{-\tau_{\rm min}}_{n+1}}. 
\end{align*}

Notice $|\widetilde f(\theta)|^{-1}$ is continuous, by Remark \ref{rmk:1-co}, we have 
\begin{align*}
\prod_{j=0}^{q_{n}-1}|\widetilde f(\theta+j\alpha)|\ge e^{q_n\big(\int_\T\ln |\widetilde f(\theta)|{\rm d}\theta-\frac{\varepsilon}{4}\big) }.
\end{align*}
By the well known integral $\int_\T\ln |\sin\pi\theta\,|{\rm d}\theta=-\ln2$, it is easy to check that 
\begin{align}
\int_\T\ln |f(\theta)|{\rm d}\theta=\int_\T\ln |\widetilde f(\theta)|{\rm d}\theta-\tau_{\rm sum}\,\ln2.
\end{align}
Therefore, 
\begin{align*}
\prod_{j=0}^{q_{n}-1}|f(\theta+j\alpha )|\ge e^{q_n\big(\int_\T\ln |\widetilde f(\theta)|{\rm d}\theta-\frac{\varepsilon}{4}\big) }e^{q_{n}\big(\delta-\frac{3\varepsilon}{4}-\tau_{\rm sum}\ln2\big)}\,{q^{-\tau_{\rm min}}_{n+1}}
 \ge e^{q_n\int_\T\ln | f(\theta)|{\rm d}\theta}\cdot \frac{e^{q_n(\delta-\varepsilon) }}{q^{\tau_{\rm min}}_{n+1}}.
\end{align*}

\qed

\subsection{Proof of Theorem \ref{thm:sc}}
Let $A(x)$ be given as in (\ref{def:A}). Direct computation shows that 
$$A^{-1}(x)=\frac{1}{f(x)} \left( \begin{array}{cc} 0 & f(x) \\ -f(x)& Ef(x)-g(x)  \end{array} \right) \triangleq \frac{F(x)}{f(x)}. $$
Observe that $F(\theta)$ and $f(\theta)$ are both  $C^{\tau_{0}}$ H\"older continuous functions. There exists cosntant $\widetilde C$ only depending the $C^{\tau_{0}}$-norms of $f$ and $g$ given in \eqref{def:norm} such that:
\begin{align}\label{eq:JY-2}
\sup_{\theta} \|F(\theta+q_{n}\alpha)-F(\theta)\|<\frac{\widetilde C}{q^{\tau_0}_{n+1}}<\frac{\widetilde C}{q^{\tau_{\rm min}}_{n+1}}, \ \ 
\sup_{\theta} \|f(\theta+q_{n}\alpha)-f(\theta)\|<\frac{\widetilde C}{q^{\tau_0}_{n+1}}<\frac{\widetilde C}{q^{\tau_{\rm min}}_{n+1}}.
\end{align}
The following lemma follows immediately from \eqref{eq:JY-1} and \eqref{eq:JY-2}. The proof is exact the same as used for Lemma 3.2 in Section 5 of \cite{JY}. We omit the details here. 
\begin{lemma}\label{lem:A}
Let $\varphi$ be a solution to $H_{\alpha, \theta}\varphi =E\varphi$ satisfying 
$\| \left(
\begin{matrix}
\varphi_0\\
\varphi_{-1}
\end{matrix}
\right)
\|=1$. 	Let $\{q_n\}$ be the subsequence given as in \eqref{eq:JY-1}, we have the following estimates:
	
	\begin{equation}\label{eq:A-1}
	\|(A^{-1}_{q_{n}}(\theta)
	-A^{-1}_{q_{n_i}}(\theta-q_{n}\alpha))
	\left(
	\begin{matrix}
	\varphi_0\\
	\varphi_{-1}
	\end{matrix}
	\right)
	\| \leq e^{q_{n}(L(E)-\delta(\alpha, \theta)+4\varepsilon)},
	\end{equation}
	and
	\begin{equation}\label{eq:A2}
	\|(A_{q_{n}}^2(\theta)-A_{2q_{n}}(\theta))
	\left(
	\begin{matrix}
	\varphi_0\\
	\varphi_{-1}
	\end{matrix}
	\right)
	\| \leq e^{q_{n}(L(E)-\delta(\alpha, \theta)+4\varepsilon)}.
	\end{equation}
	
\end{lemma}
Then combining Lemma \ref{lem:A} and the following restatement of Gordon's lemma (see e.g. \cite{JLiu,JY}), we get a contradiction, which shows the absence of point spectrum. We omit the proof here. For more details, see Sections 3 and 5 in \cite{JY}. 

\begin{thm}[Gordon's lemma]\label{thm:Gordon}
	If there exists a constant $c>0$ and a sequence $q_n$ such that the following estimates holds:
	\begin{equation}\label{absenceofpp_1}
	\|(A_{q_{n}}^2(\theta)-A_{2q_{n}}(\theta))
	\left(
	\begin{matrix}
	\varphi_0\\
	\varphi_{-1}
	\end{matrix}
	\right)
	\|\leq e^{-cq_{n}}
	\end{equation}
	and
	\begin{equation}\label{absenceofpp_2}
	\|(A^{-1}_{q_{n}}(\theta)
	-A^{-1}_{q_{n}}(\theta-q_{n}\alpha))
	\left(
	\begin{matrix}
	\varphi_0\\
	\varphi_{-1}
	\end{matrix}
	\right)
	\|\leq e^{-cq_{n}},
	\end{equation}
	then we have 
	\begin{equation}\label{maxineq}
	\max\{ \|(\begin{array}{cc}\varphi_{q_{n}} \\ \varphi_{q_{n}-1}\end{array})\|, 
	\|(\begin{array}{cc}\varphi_{-q_{n}} \\ \varphi_{-q_{n}-1}\end{array})\|,
	\|(\begin{array}{cc}\varphi_{2q_{n}} \\ \varphi_{2q_{n}-1}\end{array})\| \} \geq \frac{1}{4}.
	\end{equation}
\end{thm}

\section{Generic full spectral/packing dimension}\label{sec:dim}
In this part, we fix $\theta$ such that $\delta(\alpha,\theta)= \tau_{\rm min}\cdot\beta(\alpha)$ as in \eqref{def:delta}. Let the subsequence $q_n$ the given as in \eqref{eq:q_n} and \eqref{eq:JY-1} satisfying
\begin{align}\label{eq:q_n2}
\delta(\alpha,\theta)>\frac{\tau_{\rm min}\ln q_{n+1}} {q_{n}} >\delta(\alpha,\theta)-\frac{\varepsilon}{4}\Longleftrightarrow e^{q_n \tau_{\rm min}\cdot\beta }>q^{\tau_{\rm min}}_{n+1}>e^{q_n(\tau_{\rm min}\cdot\beta-\varepsilon/4)}. 
\end{align}
By repeating the argument in the proof of Lemma \ref{lem:key} and \ref{lem:key2}, one can prove a more general version of (\ref{eq:JY-2}) for $\theta+k\alpha$, $|k|<e^{\varepsilon q_n}$. The same type of estimates has been obtained in the recent preprint of \cite{HYZ} for a.e. $\theta$ such that $\delta(\alpha,\theta)= \tau_{\rm min}\cdot\beta(\alpha)$. We restate the result here directly. Readers can find more details in Lemma 4.2 in \cite{HYZ}.  
\begin{lem}
Let $f\in{\cal F}$ be given as in \eqref{def:F} with $\int_\T\ln | f(\theta)|{\rm d}\theta=0$. For $\beta(\alpha),\varepsilon>0$ , let the sequence $q_n\to \infty$  be defined as in \eqref{eq:q_n2}. There is a full Lebesgue measure set $\Theta=\Theta(\alpha,\theta_1,\cdots,\theta_m)$ such that 
	for any $\theta\in\Theta$ and $q_n$ large enough\footnote{The sequence itself only depends on $\beta(\alpha)$, while the largeness depends on $\theta,\alpha,\beta,\varepsilon,\tau$.},  $f(\theta+n\alpha)$ satisfies:
	\begin{align}\label{eq:q_n3}
	\min_{|m|\le e^{\varepsilon q_n/10}}\prod_{j=m}^{m+q_n-1}   |f(\theta+j\alpha)|>e^{-\varepsilon\,q_n}.
	\end{align}
\end{lem}

Let $D$ and $D_n$ be defined as in \eqref{def:Mn},\eqref{def:AD}.  By Lemma \ref{lem:F}, for $r$ large enough, 
\begin{align}\label{eq:Df}
\|D_r(\theta)\| \leq e^{(L(D)+\varepsilon) r}.  
\end{align}

Combine the definition (\ref{def:AD}) and the estimates \eqref{eq:q_n3},\eqref{eq:Df}, direct computation shows the following Lyapunov growth of the unbounded cocyle $(\alpha,A)$. For $E$ in a compact set $I$, let $\Lambda=2\varepsilon+\sup_{E\in I}L(E)$, we have
\begin{lem}\label{lem:infA}
	For any $\theta\in\Theta$ and $q_n$ large enough and $1\le r\le q_n$, 
	\begin{align}\label{eq:Au}
\sup_{{|m|\le e^{\varepsilon q_n/10}}}	\|A_r(\theta+m\alpha)\|< e^{(L(E)+2\varepsilon) q_n}=e^{\Lambda q_n}. 
	\end{align}
\end{lem}
\begin{rmk}
	The above lemma actually shows that the following $\liminf$
	\begin{align}\label{eq:liminf}
	\liminf_n\, \frac{1}{n}\ln\|A_n(E,\theta)\|\le L(E)
	\end{align}
	holds uniformly on a full measure set 
	$\widetilde \Theta=\bigcup_{|m|\le e^{\varepsilon q_n/10}}(\Theta+m\alpha)$. It was showed in \cite{F} that the $\limsup$ in \eqref{eq:liminf} holds uniformly on $\T$ for any continuous cocycle. The uppersemi continuity was generalized in \cite{JM} to almost continuous bounded cocyles. In general, the $\limsup$ in \eqref{eq:liminf} shall not hold for unbounded cocyles. Lemma \ref{lem:infA} can be viewed as a quantitative generalization of  \cite{F,JM} to the unbounded case. 
\end{rmk}
 
On the other hand, by (\ref{eq:JY-2}) and \eqref{eq:q_n3},  $A_n(\theta)$ has strong repetitions ($\tau_{\rm min}\beta$-almost periodicity): \eqref{eq:JY-2}:
\begin{align}\label{eq:Abeta}
\sup_{\theta} \|A(\theta+q_{n}\alpha)-A(\theta)\|<\frac{\widetilde C}{q^{\tau_{\rm min}}_{n+1}}<\widetilde Ce^{-q_n(\tau_{\rm min}\beta-\varepsilon/4)}. 
\end{align}
As proved in \cite{JZ}, the growth of the transfer matricies \eqref{eq:Au} and the strong repetition \eqref{eq:Abeta} imply that there exists an absolute constant $C_1$ such that $\mu_I$ is $\gamma$-spectral continuous for all $\gamma\le 1-\frac{C_1\Lambda}{\tau_{\rm min}\beta-\varepsilon/4}$. This completes the proof of Theorem \ref{thm:dim} by the definition of spectral dimension in \eqref{def:dimspe}. We omit the proof here. For more details about the $\gamma$-spectral continuity, we refer readers to Theorem 1 and Theorem 6 in \cite{JZ}. 

 \section*{Acknowledgement}
 The authors would like to thank Rui Han for useful discussions. The authors would also like to thank Ilya Kachkovskiy for mentioning important literatures of unbounded Schr\"odinger operators to us. Last but not least, the authors would like to thank Svetlana Jitomirskaya for reading the early manuscript and useful 
comments. F. Y. would like to thank the Institute for Advanced Study, Princeton, for its hospitality
 during the 2017-18 academic year.  F. Y. was supported in part by NSF grant DMS-1638352.
 S. Z. was supported in part by NSF grant DMS-1600065.

\bibliographystyle{amsplain}

\vspace{1cm}
Fan Yang,

School of Math, Institute for Advanced Study.

E-mail address: yangf@ias.edu

\vspace{1cm}
Shiwen Zhang, 

Dept. of Math., Michigan State University. 

E-mail address: zhangshiwen@math.msu.edu

\end{document}